\magnification=1000
\documentstyle{amsppt}
\parindent 20 pt
\baselineskip 20pt
\define \bul{\bullet}

\define \hK {$\hat K$}

\define \a{\alpha}

\define \r{\rho}
\define \s{\sigma}

\define \z{\zeta}

\define \iy{\infty}

\define \ri{\rightarrow}

\define \sbt{\subset}
\define \spt{\supset}

\define \edm{\enddemo}
\define \ep{\endproclaim}

\define \1{^{-1}}
\define \2{^{-2}}

\define \BC{\Bbb C}

\define \BR{\Bbb R}

\define \Cd{\Cal D}

\define \var{\operatorname{Var\ }}

\define \intl{\int\limits}

\def\eps{\epsilon}

\define \Fr {\bold {Fr}}
\define \Di{\Cd}
\topmatter \heading{\bf Limit sets and a problem in dynamical systems   }\endheading
 \vskip .10in
\centerline{Vladimir Azarin }
\NoBlackBoxes
\abstract{We approximate  a chain recurrent dynamical system by periodic
dynamical systems. This is similar to the well known Bohr theorem on
approximation of almost periodic functions by periodic functions. }
\endabstract
\endtopmatter

\heading{1.Introduction and the main results}\endheading
\subheading{1.1. Chain recurrence}

A family of maps of the form
$$ T^t:M\mapsto M,\ t\in \BR$$
on a compact connected metric space $M$ with a metric $d(\bullet,\bullet)$ is a {\it dynamical system} $(T^\bullet,M)$ if it satisfies the condition
$$T^{t+\tau}=T^t\circ T^\tau, \forall t,\tau\in \BR$$
and the map $(t,m)\mapsto T^tm $ is continuous with respect to $(t,m)\
\forall t\in \BR,m\in M$ in the natural topology.

Let $m,m'\in M,$ the number $\eps >0$ be small and the number $s>0$
be large. An $(\eps,s)$-chain from $m$ to $m'$ is a finite sequence
$m_0=m, m_1,...m_n=m'$ satisfying the conditions
$d(T^{t_j}m_j,m_{j+1})<\eps,$ for some $t_j>s ,\ j=0,1,...,n-1.$

A dynamical system $(T^\bullet, M)$ is called {\it chain recurrent}, if for arbitrary small $\eps>0$ and arbitrary large $s>0$ there exists an
$(\eps,s)$-chain in $M$ from $m$ to $m.$

\subheading{1.2. ADPT}

Let $m(t):[0,\iy) \mapsto M$ be a curve (pseudo-trajetory) in $M.$ It is called
{\it dense} in $M$ if
$$clos \{m(t):t\in [a,\iy)\}=M, \forall a\in \BR$$
where $clos $ means closure.

A pseudo-trajectory is called {\it asymptotically dynamical with the dynamical
asymtotics} $T^\bullet$ if
$$ d(T^{t+\tau}m(t),m(t+\tau))\ri 0$$
as $t\ri \iy$ uniformly with respect to $\tau\in [a,b],\forall a,b<\iy.$

The curve can be piecewise continuous.

The following assertion holds:
\proclaim {Theorem 1.1} $(T^t,M)$ is chain recurrent iff there exists a dense a.d.p.t.\ep

We omit its proof (see e.g.\cite {Az, Th.4.3.1.2, 4.3.3.3.}).
\subheading {1.3.Approximation}
Let $M$ be a metric compact space with a   distance $d(\bullet,\bullet).$ For a set $ E\sbt M$ define
$$E_\eps:=\{m\in M:\exists e\in E :d(m,e)<\eps\}$$
-- a neighborhood of $E.$

For sets $E',E''\sbt M$ define
$$d(E,E'):=\inf \{\eps: E'\sbt E''_\eps \wedge E''\sbt E'_\eps\}.$$
We write
$E_n\ri E$ if $d(E_n,E)\ri 0.$

Let $M'\spt M$ be a compact set and let $(T^{\bullet '} ,M')$ be a dynamical system.
We say that $(T^{\bullet '},M')$ {\it reduces} to $M$ as $(T^\bullet,M)$ if
$(T^{\bullet '},M)=(T^\bullet,M).$
\proclaim {Theorem 1.2. (Approximation)} Let $(M,T^t),t\in (-\iy,\iy)$ be a chain recurrent dynamical system on
a connected metric compact set $M$, and let $$Orb(x):=\{T^tx: -\iy\leq t\leq \iy\}$$
be the orbit passing through a point $x\in M.$

Suppose the orbit is not periodic.
Then there exists a dynamical system $(T^{\bullet '},M')$, $M'\spt M$ that reduces to
$M$ as $(T^t,M)$, and for every sequence of numbers $P_n\ri\iy$ there exists a sequence
of periodic orbits
$$Orb(x_n):=\{T^{t'}x_n: 0\leq t\leq 2P_n\},\ x_n\ri x$$
such that $Orb(x_n)\ri \ Orb(x).$
\ep

\heading 2.Proofs\endheading

\subheading {2.1. Model}

Now we realize the chain recurrent dynamical system by the following model.

Consider the set $\Cal M[\r,\s]$ of positive measures $\mu(E),E\sbt \BC,$ satisfying the condition:

$$\Cal M[\r,\s]:=\{\mu: \mu(r):=\mu(|\z|<r)\leq \s r^\r \forall r>0\}$$
It is metrizable as a compact set in $\Di' $ -topology, over the space
$\Di (\BC\setminus 0)$ of infinitely differentiable functions that are finite in $\BC\setminus 0.$
This topology is metrizable ,because it is a {\it Fr\'echet} space.
 Denote the distance as $d(\bullet,\bullet).$

$\Cal M[\r,\s],\forall \r>0,\s>0$  is invariant with respect to the transformation
$$T_t:T_t\mu(E)=\mu (e^t E)e^{-\r t}.$$
$T_t m$ is continuous in $(t,\mu), \ t\in \BR,\mu\in \Cal M[\r,\s]$ in the product of
the corresponding topologies.
So the dynamical system $(T_\bullet,\Cal M[\r,\s])$ is defined.

Let $\mu\in \Cal M[\r,\s].$ Truncate it: let $\mu^*_{P_n}$ be the reduction
to the annulus $e^{-P_n}\leq |z|\leq e^{P_n},$; now extend it periodically to all
$\BC$ in the following way:
 $$\mu_{P_n}:=\sum\limits_{k=-\iy}^{\iy}T_{k2P_n}\mu^*_{P_n}.$$
Every orbit passing through $ \mu_{P_n}$ is periodic with period $2P_n.$

\proclaim {Proposition 2.1}The sequence $\mu_{P_n}\ri \mu$ in $ D'$ topology\ep and, hence,
in the equivalent metric.

\demo {Proof}This is so because for every $g\in \Di(\BC\setminus 0)$ we have
$supp \, g\sbt (e^{-P_n},e^{P_n})$ for large $n.$ Hence
$$<\mu_{P_n},g>:=\int g(z)\mu_{P_n}(dz) = \int g(z)\mu(dz)$$
for large $n.$
\edm
 \proclaim{Theorem 2.2 (Approximation)}
$$Orb(\mu_{P_n})\ri Orb(\mu).$$
\ep
\demo {Proof} In the definition of $d(X_n,X)$ (see \S1.3) set
$X_n:=Orb(\mu_{P_n}),X=Orb(\mu).$ If $d(X_n,X) \not\ri 0,$ we get a contradiction
to Theorem 2.1.
\edm
So the Theorem 1.3  is proved for $M=\Cal M[\r,\s]$ and $T^t=T_t.$

Now we should prove
\proclaim {Theorem 2.3 (Universality of $\Cal M[\r,\s]$)}Let $(T^\bul,M)$ be
a chain recurrent dynamical system on a compact set $M.$ Then for any $\r,\s$ there exists
$\Cal M\sbt \Cal M[\r,\s]$ and a homeomorphism $imb:M\mapsto \Cal M$ such that
$imb\circ T^t=T_t\circ imb,\ t\in (-\iy,\iy).$
\ep
i.e., any dynamical system can be imbedded in $(T_\bul,\Cal M[\r,\s]).$

This  theorem is proved in \cite {Az ,Th.4.1.5.1}

Here is a sketch of the proof.

Let us denote by $\Cal M (S)$ the set of measures $\nu$ with bounded full variation on the
unit circle $S.$ Introduce the metric $d(\nu,0):=\var \nu$ and consider the set
$$K:=\{\nu:\nu>0,\ d(\nu,0)\leq 1\},$$
i.e., the intersection of the cone of positive measures with the unit ball.

The following assertion is a corollary of the Keller's theorem  (see,e.g.
\cite {BP, Th.3.1, p.100}).
\proclaim {Theorem 3.3 (Imbedding)} Every metric compact set can be
homeomorphically imbedded in $K.$
\ep
Thus we can assume below that for any $m\in M$ there exists a positive
measure
$$ Y(\bul,m)=Y(d\phi,m)\in K $$
\demo {Proof of Universality Theorem}The numeration corresponds to that of
\cite {Az,\S 4.1.6}

Let us ``transplant'' $\mu $ in
$$Cyl:=S\times\BR~.$$ For $\mu \in\Cal M[\r,\s]$ that has a density $f_\mu (re^{i\phi}),$
we set
$$\nu (dy\otimes d\phi):=f_\mu (e^ye^{i\phi})e^{-\r y}(dy\otimes d\phi).$$
i.e., the density $f_\nu$ of $\nu$ is defined by
$$f_\nu(\phi,y):=f_\mu (e^ye^{i\phi})e^{-\r y}.$$
Accordingly,
$$f_\mu (\phi,r)=f_\nu(\phi,\log r)r^{\r}$$
We can extend this equality to all $\mu\in \Cal M[\r,\s]$ using a
limit process in $\Di '$ topology.

We can also define $\nu$ as a distribution in $\Di'(Cyl).$ Namely, for
$\psi\in \Di(Cyl)$ we set
$$\psi^*(\phi,r):=\psi (Pol^{-1}(\phi,\log r))r^{-\r y}~,$$
where :
$$Pol:\BC\setminus 0\mapsto (0,\iy)\times S_1$$
by the polar coordinates $re^{i\phi}\mapsto (r,\phi).$
Let $X(t)$ be a positive function  satisfying the condition
$$\intl_{-\iy}^{\iy}X(t)dt=1$$
and such that the linear hull of its
translations are dense in $L^1(-\iy,\iy)$.
 We can choose, for example, the
function
$$X(t):=\frac {1}{\sqrt {2\pi}}e^{-\frac{t^2}{2}}$$
because its Fourier transformation does not vanish in $\BR$
(it is $e^{-\frac{ s^2}{2}}$).

Let us define $\nu (\bul,m)$ by
$$<\nu(\bul,m) ,\psi):=\intl_{(\phi,y)\in Cyl} \psi (\phi,y)\left (\r\intl_{-\iy}^{\iy}Y(d\phi,T^{y-t}m)X(t)dt\right ) dy.\tag 4.1.6.4$$
 Now we check the property
$$S_\tau \nu(\bul,m)=\nu (\bul, T^\tau m)$$
Using (4.1.6.2), we obtain
$$<S_\tau\nu (\bul, m),\psi>=
\int \psi (\phi,y)
\left (\r\intl_{-\iy}^{\iy}Y(d\phi,T^{y+\tau-t}m)X(t)dt\right ) dy=$$
$$\int \psi (x^0,y)
\left (\r\intl_{-\iy}^{\iy}Y(d\phi,T^{y-t}( T^\tau m))X(t)dt\right ) dy=
<\nu (\bul, m),T^\tau m>.$$

and
$$<\nu,\psi>:=\int \psi^* (\phi, r)\mu (d\phi\otimes rdr)$$
The transformation $P_t x=(\phi,tr),\ re^{i\phi}\in \BC\setminus 0$ passes to
$$Pol\circ P_t\circ Pol^{-1}(\phi,y)=( \phi,y+\log t)$$
Thus $T_t\mu$ gives a transformation $S_t \nu$ defined by
$$S_t f_\nu (\phi,y):=f_\nu (\phi,y+t)$$
for densities or by
$$<S_t\nu,\psi>:=\int \psi (\phi,y-t)\nu (d\phi\otimes dy)\tag 4.1.6.2$$
for distributions ( $\psi \in \Di (Cyl).$)

From $\mu\in \Cal M[\r,\s]$  we obtain
$$\intl_{y\leq 0}e^{\r y}S_t\nu (dy\otimes d\phi)\leq \s,\  t\in\BR~,$$

So we should check the one-to-one correspondence between $\nu(\bul,m)$ and $Y(\bul,m).$

Suppose
$$\nu(\bul,m_1)=\nu(\bul,m_2).$$
Then
$$<\nu(\bul,m_1),\psi>=<\nu(\bul,m_2),\psi>\ \forall \psi\in \Di(Cyl).$$
In particular, set
$$\psi(\phi,y)=\Phi (\phi)R(y),\ \Phi\in \Di (S),\ R\in \Di (-\iy,\iy).$$
Then
$$<\nu(\bul,m_1),\psi>=\int R(y)dy
\intl_{-\iy}^\iy <Y(\bul,T^{y-t}m_1),\phi>_{S}
X(t)dt=\tag 4.1.6.6$$
$$=<\nu(\bul,m_2),\psi>=\int R(y)dy
\intl_{-\iy}^\iy <Y(\bul,T^{y-t}m_2),\phi>_{S}.$$
where
 $$ <Y(\bul),\phi>_{S}:=\intl_{S}\phi(\phi)Y(d\phi).$$
Set
$$F_j(y):=<Y(\bul,T^y m_j),\phi>_{S},\ j=1,2.$$
From  (4.1.6.6) we obtain for the convolutions
$$(F_1*X)(y)\equiv(F_2*X)(y),\ y\in (-\iy,\iy).$$
Thus
$$F_1(y)\equiv F_2(y),\ y\in (-\iy,\iy)$$
because of the property of $X.$

Hence
$$Y(\bul,T^y m_1)\equiv Y(\bul,T^y m_2),\ y\in (-\iy,\iy).$$
In particular, for $y=0$ we have
$$Y(\bul,m_1)=Y(\bul,m_2).$$
Hence $m_1=m_2$ because of (4.1.5.3), and this completes the proof of one-to-one
correspondence.
\qed
\edm
So we proved also Theorem 1.2 (Approximation).

Let us note that the approximating  dynamical system is considered
on a larger compact set than $Orb(x),$ but reducing it to the
closure of $Orb(x)$ yields the same dynamical system. \subheading
{2.3.Example} Consider an example that was proposed to us as a
counterexample by Prof.Morris W.Hirsh and was a reason to change a
little the formulation.

The  first realization:Let $M=S^1\times S^1$ and $(\phi_x,\theta_x)$ define every point in $M.$ Set
$T^tx=(\phi_x +2\pi t,\theta_x+\a 2\pi t)$ where $\a$ is an irrational number.

This DS is homeomorphic to the following one.

The second realization:place in  points $x_1,x_2$
two masses $m_1=0.5,m_2=0.5$, this is a mass distribution $\mu _0$ on $S$
($\mu _0\in K$) and consider as $T^t\mu_0$ moving this mass distribution as described above.

Now consider $Hom:K\mapsto \Cal M[\r,\s]$ such that $\mu _0\mapsto\mu$ which
is the mass distribution on the positive  ray that consists of two masses
$m_1=(1/2)\s;m_2=(1/2)\s-\eps$ concentrated , for example, at points
 $x=1;y=0,$ and $x=\a,y=0.$ Then $\mu_{P_n}\in \Cal M[\r,\s]$ and
the assertion of Approximation Theorem holds.

Returning to the second realization we obtain that the irrational curve on the torus can be approximated
by an irrational curve (with masses) on an infinite-dimensional torus.

I am indebted to Prof.A.Eremenko who pointed out to me a mistake in the previous formulation.

\end